\documentclass{amsart}
\usepackage{amssymb}
\usepackage{amscd}
\usepackage{amsthm}

\newtheorem{theo}{Theorem}[section]

\newtheorem{prop}[theo]{Proposition}
\newtheorem{lem}[theo]{Lemma}
\newtheorem{cor}[theo]{Corollary}

\newtheorem{defi}[theo]{Definition}

\theoremstyle{remark}
\newtheorem{rema}[theo]{Remark}
\newtheorem{sect}[theo]{}

\def \Romannumeral #1 {\expandafter\uppercase\expandafter
{\romannumeral #1} }

\def \Br {{\rm{Br\,}}}

\def \Ga {{\Gamma}}
\def \R {{\bf R}}
\def \Pic {{\rm {Pic\,}}}

\def \Gal {{\rm{Gal\,}}}

\def \A{{\mathbb A}}

\def \dim {{\rm{dim\,}}}
\def \Hom {{\rm {Hom}}}
\def \End {{\rm {End}}}

\def \Aut{{\rm Aut \,}}

\def\ov{\overline}

\def \Z {{\bf Z}}
\def \Q {{\bf Q}}
\def \F {{\bf F}}
\def \Br {{\rm Br\,}}

\def \sym {{\rm{sym}}}

\def \rk {{\rm{rk}}}
\def \non {{\rm{non-}}}
\def\G{{\bf G}}
\def\R{{\bf R}}

\def\C{{\bf C}}

\def\lra{\longrightarrow}

\def\H{{\rm H}}

\def\NS{{\rm NS\,}}

\def\Ga{\Gamma}

\def\et{{\acute et}}
\def\tors{{\rm tors}}

\def\fchar{{\rm char}}
\def\Lie{\mathrm{Lie}}
\newcommand{\bthe}{\begin{theo}}
\newcommand{\ble}{\begin{lem}}
\newcommand{\bpr}{\begin{prop}}
\newcommand{\bco}{\begin{cor}}
\newcommand{\bde}{\begin{defi}}
\newcommand{\ethe}{\end{theo}}
\newcommand{\ele}{\end{lem}}
\newcommand{\epr}{\end{prop}}
\newcommand{\eco}{\end{cor}}
\newcommand{\ede}{\end{defi}}

\def\bsm{\left( \begin{smallmatrix}}
\def\esm{\end{smallmatrix} \right)}

\title[Brauer groups of abelian
varieties and K3 surfaces]{A finiteness theorem for the Brauer
group of abelian varieties and K3 surfaces}
\author[Alexei\ N.\ Skorobogatov]{Alexei\ N.\ Skorobogatov}
\address{Department of Mathematics, South Kensington Campus,
Imperial College London, SW7 2BZ England, U.K.}
\address{Institute for the Information Transmission Problems, Russian
Academy of Sciences, 19 Bolshoi Karetnyi, Moscow, 127994 Russia}
\email{a.skorobogatov\char`\@imperial.ac.uk}
\author[Yuri\ G.\ Zarhin]{Yuri\ G.\ Zarhin}
\address{Department of Mathematics, Pennsylvania State University,
University Park, PA 16802, USA}
\address{Institute for Mathematical Problems in Biology,
Russian Academy of Sciences, Pushchino, Moscow Region, Russia}
\email{zarhin\char`\@math.psu.edu}

\begin{document}

\baselineskip=15pt
\begin{abstract}
Let $k$ be a field finitely generated over the field of
rational numbers, and
$\Br(k)$ the Brauer group of $k$. For an algebraic
variety $X$ over $k$ we
consider the cohomological Brauer--Grothendieck group $\Br(X)$.
We prove that the quotient of $\Br(X)$ by the image of $\Br(k)$
is finite if $X$ is a $K3$ surface. When $X$ is
an abelian variety over $k$,
and $\ov X$ is the variety over an algebraic closure $\ov k$ of $k$ 
obtained from $X$ by the extension of the ground field,
we prove that the image of $\Br(X)$ in $\Br(\ov X)$ is finite.
\end{abstract}

 \maketitle

\section{Introduction}

Let $X$ be a geometrically integral smooth projective variety
over a field $k$.
The Tate conjecture for divisors on $X$ \cite{tateW,tate,T91} is well
known to be closely related to the finiteness properties of the
cohomological Brauer--Grothendieck group
$\Br(X)=\H^2_\et(X,\G_m)$. This fact was first discovered 
in the case of a finite field $k$
by Artin and Tate (\cite{tateB}, see
also Milne \cite{Milne}) who studied the Brauer group of a
surface. In particular, the order of $\Br (X)$ appears in the
formula for the leading term of the zeta-function of $X$.
A stronger variant of
the Tate conjecture for divisors concerns the order of the
pole of the zeta-function of $X$ at $s=1$, see \cite[(12) on p. 101]{tateW}.
It implies the finiteness of the prime-to-$p$ component of $\Br(X)$,
where $X$ is a variety of arbitrary dimension,
and $k$ is a finite field of characteristic $p$, as
proved in \cite[Sect. 2.1.2 and Remark 2.3.11]{Z82}.

Since Manin observed that the Brauer group of a variety over a
number field provides an obstruction to the Hasse principle
\cite{Manin}, the Brauer groups of varieties over fields of
characteristic $0$ have been intensively studied. Most of the
existing literature
is devoted to the so called algebraic part $\Br_1(X)$ of
$\Br(X)$, defined as the kernel of the natural map
$\Br(X)\to\Br(\ov X)$, where $\ov X=X\times_k\ov k$, and $\ov k$
is a separable closure of $k$.
Meanwhile, if $k$ is a number field, the
classes surviving in $\Br (\ov X)$ can produce a non-trivial
obstruction to the Hasse principle and weak approximation (see
\cite{H} and \cite{W} for explicit examples). Therefore,
such arithmetic applications require the knowledge of
the whole Brauer group $\Br (X)$.

To state and discuss our results we introduce 
some notation and conventions. 
In this paper the expression `almost all' means `all but finitely many'.
If $B$ is an abelian group, we denote by $B_\tors$ the torsion
subgroup of $B$, and write $B/_\tors:=B/B_\tors$.
For a prime $\ell$ let $B(\ell)$ be the subgroup of $B_\tors$
consisting of the elements whose order is a power of
$\ell$, and $B({\non}\ell)$ the subgroup of $B_\tors$ consisting
of the elements whose order is {\sl not} divisible by $\ell$. If $m$
is a positive integer, we write $B_m$ for the kernel of the
multiplication by $m$ in $B$.

Let $\Br_0(X)$ be the image of the natural map $\Br(k)\to\Br(X)$.
Recall that both $\Br (X)$ and
$\Br(\ov X)$ are torsion abelian groups whenever $X$ is smooth, see
\cite[II, Prop. 1.4]{Gr}. There are at least three reasons
why the Brauer group $\Br (X)$ can be infinite:
$\Br_0 (X)$ may well be infinite; 
the quotient $\Br_1(X)/\Br_0(X)$ injects into, 
and is often equal to,
$\H^1(k,\Pic(\ov X))$, which may be infinite if the divisible part
of $\Pic(\ov X)$ is non-zero,
or if there is torsion in the N\'eron--Severi group $\NS(\ov X)$;
finally, $\Br(\ov X)$ may be infinite.
This prompts the following question.

{\it Question $1$. Is $\Br(X)/\Br_1(X)$ finite if $k$ is finitely
generated over its prime subfield?}

\noindent Let $\Ga=\Gal(\ov k/k)$, and 
let $\Br(\ov X)^{\Ga}$ be the subgroup of
Galois invariants of $\Br(\ov X)$; then $\Br (X)/\Br_1 (X)$ naturally
embeds into $\Br(\ov X)^{\Ga}$. A positive answer to Question 1
would follow from a positive answer to the following question.

{\it Question $2$. Is $\Br(\ov X)^{\Ga}$ finite if $k$ is finitely
generated over its prime subfield?}

\noindent
In this note we prove the following two theorems.

\bthe \label{t1} Let $k$ be a field finitely generated over its prime
subfield. Let $X$ be a principal homogeneous space of an abelian
variety over $k$.
\begin{itemize}
\item[(i)] If the characteristic of $k$ is $0$, then the groups
$\Br(\ov X)^{\Ga}$ and $\Br (X)/\Br_1 (X)$ are finite.

\item[(ii)] If the characteristic of $k$ is a prime
$p\not=2$, then the groups
$\Br(\ov X)^{\Ga}({\non}p)$ and $(\Br (X)/\Br_1 (X))({\non}p)$ are finite.
\end{itemize}
\ethe

\bthe \label{t} Let $k$ be a field finitely generated over $\Q$.
If $X$ is a K3 surface over $k$, then the groups 
$\Br(\ov X)^{\Ga}$ and $\Br (X)/\Br_0 (X)$ are finite. 
\ethe
\begin{rema}
\label{GammaIn} The injective maps 
$$\Br (X)/\Br_1 (X)\hookrightarrow
\Br(\ov X)^\Ga\quad \text{and}\quad \Br_1(X)/\Br_0(X)\hookrightarrow
\H^1(k,\Pic(\ov X))$$ 
can be computed via the Hochschild--Serre spectral sequence
$$
\H^p(k,\H_\et^q(\ov X,\G_m))\Rightarrow \H_\et^{p+q}(X,\G_m).
$$
(A description of some of its differentials can be found in
\cite{S1}.) 
Recall that in characteristic zero the Picard group
$\Pic(\ov X)$ of a K3
surface $X$ is a free abelian group of rank at most
$20$. The Galois group $\Ga$ acts on
$\Pic(\ov X)$ via a finite quotient, so that
$\H^1(k,\Pic (\ov X))$ is finite. Thus in order to prove Theorem
\ref{t} it suffices to establish the finiteness of
$\Br(\ov X)^\Ga$.
\end{rema}

In the case when the rank of $\Pic(\ov X)$ equals $20$, Theorem
\ref{t} was proved by Raskind and Scharaschkin \cite{RS}. In an
unpublished note J-L. Colliot-Th\'el\`ene proved that $\Br(\ov
X)^{\Ga}(\ell)$ is finite for every prime $\ell$, where $X$ is a
smooth projective variety over a field finitely generated over
$\Q$, assuming the Tate conjecture for divisors on
$X$. (When $\dim(X)>2$, he assumed additionally the semisimplicity
of the Galois action on $\H_\et^2(\ov X,\Q_\ell)$.)
See also \cite{Tankeev} for some related results.

When $X$ is an abelian variety over a field finitely generated
over its prime subfield, the Tate conjecture for
divisors on $X$ (and the semisimplicity of $\H_\et^2(\ov X,\Q_\ell)$
for $\ell \ne p$) was proved by the second named author
in characteristic $p>2$
\cite{Z75,Z76}, and by Faltings in characteristic zero
\cite{F,F86}. This result of Faltings combined with the construction of
Kuga--Satake  elaborated by Deligne \cite{D}, implies the Tate
conjecture for divisors on K3 surfaces in characteristic zero
\cite[p. 80]{T91}.

The novelty of our approach is due to the usage of a variant of
the Tate conjecture for divisors on $X$ \cite{Z77,Z85} which concerns
the Galois invariants of the (twisted) second \'etale cohomology
group with coefficients in $\Z/\ell$ (instead of $\Q_{\ell}$),
for almost all primes $\ell$. Using this variant we
prove that under the conditions of Theorems \ref{t1} and \ref{t}
we have
$\Br(\ov X)^{\Ga}_{\ell}=\{0\}$ for almost all primes $\ell$.

\medskip

Let $k$ be a number field, $X({\mathbb{A}}_k)$
the space of adelic points of $X$, and $X({\mathbb{A}}_k)^\Br$
the subset of adelic points orthogonal to $\Br(X)$ with respect to the
Brauer--Manin pairing (given by the sum of local invariants
of an element of $\Br(X)$ evaluated at the local points, see \cite{Manin}).
We point out the following corollary to Theorem \ref{t}.

\bco \label{last} Let $X$ be a K3 surface over a number field $k$.
Then $X({\mathbb{A}}_k)^\Br$ is an open subset of
$X({\mathbb{A}}_k)$. \eco
\begin{proof} The sum of local invariants of a given
element of $\Br(X)$ is a continuous function on
$X({\mathbb{A}}_k)$ with finitely many values. 
Thus the corollary is a consequence of
Theorem \ref{t}. \end{proof}

Let us mention here
some open problems regarding rational points on K3 surfaces.
Previous work on surfaces fibred into curves of genus 1
\cite{CSS,SD,SS} indicates that it is not unreasonable to expect
the Manin obstruction to be the only obstruction to the Hasse
principle on K3 surfaces. One could raise a more daring question:
is the set of $k$-points dense in the Brauer--Manin set
$X(\A_k)^\Br$? By Corollary \ref{last}, this would imply that the
set of $k$-points on any K3 surface over a number field is either
empty or Zariski dense. Moreover, this would also imply the weak-weak
approximation for $X(k)$, whenever this set is non-empty. (This means that
$k$ has a finite set of places $S$ such that for
any finite set of places $T$ disjoint from $S$
the diagonal image of $X(k)$ in $\prod_{T}X(k_v)$ is dense.)

\medskip

The paper is organized as follows. In Section \ref{prel} we recall the
basic facts about the interrelations between the Brauer group, the
Picard group and the
N\'eron--Severi group (mostly due to Grothendieck \cite{Gr}).
We also discuss some linear algebra constructions arising from
$\ell$-adic cohomology. In Section \ref{TateC} we recall the
finite coefficients variant of the Tate conjecture for abelian varieties,
and prove Theorem \ref{t1}. 
Finally, Theorem \ref{t} is proved in Section \ref{K33}.

We are grateful to the referee whose comments helped to improve
the exposition. This work was done during the special semester
``Rational and integral points on higher-dimensional varieties" at
the MSRI, and we would like to thank the MSRI and the organizers
of the program. The second named author is grateful to
Dr. Boris Veytsman for his help with \TeX nical problems.

\section{The N\'eron--Severi group, $\H^2$ and the Brauer group}
\label{prel}

We start with an easy lemma from linear algebra.

\begin{lem}
\label{discr} Let $\Lambda$ be a principal ideal domain, $H$ a
non-zero $\Lambda$-module, $N\subset H$ a non-zero free submodule of finite
rank. Let
$$\psi: H \times H \to \Lambda$$
be a symmetric bilinear form. Let $N^\perp$ be the orthogonal
complement to $N$ in $H$ with respect to $\psi$, and let $\delta$ be the
discriminant of the restriction of $\psi$ to $N$.
If $\delta\not=0$, then $N\cap N^\perp=\{0\}$ and 
$$\delta^2 H\subset N\oplus N^\perp\subset H.$$ 
In particular, if $\delta$ is a
unit in $\Lambda$, then $H=N\oplus N^\perp$.
\end{lem}
\begin{proof}
Let us put $N^{*}=\Hom_{\Lambda}(N,\Lambda)$. The form $\psi$
gives rise to a natural homomorphism of $\Lambda$-modules
$e_{\psi}:H \to N^*$ with $N^\perp=\ker(e_{\psi})$ and
$$\delta\cdot N^{*}\subset e_{\psi}(N)\subset N^{*}.$$
In particular, the restriction of $e_{\psi}$ to $N$ is injective,
therefore $N\cap N^\perp=\{0\}$, and  $e_{\psi}:N\to e_{\psi}(N)$
is an isomorphism. Let $u:e_{\psi}(N)\cong N$ be its inverse,
i.e., $u e_{\psi}:N\to N$ is the identity map. Let us consider the
homomorphism of $\Lambda$-modules
$$P:H \to N, \quad h \mapsto \delta u (e_{\psi}(h)).$$
This definition makes sense
since $\delta e_{\psi}(h)\in \delta N^{*}\subset e_{\psi}(N)$.
It is clear that $\delta \cdot \ker{P}\subset N^\perp\subset \ker(P)$,
and $P$ acts on $N$ as the
multiplication by $\delta$. For any $h \in H$ we have $z=P(x)\in N$
and $P(z)=\delta z$, which implies that $P(\delta h)= P(z)$. Hence
$\delta h -z \in \ker(P)$, and therefore $\delta(\delta h -z) \in
N^\perp$. It follows that $\delta^2 h \in \delta z+N^\perp\subset N\oplus
N^\perp$.
\end{proof}

\begin{sect}
\label{kummer}
Let us recall some useful elementary statements,
which are due to Tate \cite{tateB,Tate76}.
Let $B$ be an abelian group. The
projective limit of the groups $B_{\ell^n}$ (where the transition maps
are the multiplications by $\ell$) is called the $\ell$-adic 
Tate module of $B$, and is denoted by $T_{\ell}(B)$. This
limit carries a natural structure of a  $\Z_{\ell}$-module;
there is a natural injective map $T_{\ell}(B)/\ell
\hookrightarrow B_{\ell}$. One may easily check that
$T_{\ell}(B)_{\ell}=\{0\}$, and therefore $T_{\ell}(B)$ is
torsion-free. Let us assume that $B_{\ell}$ is finite. Then all the
$B_{\ell^n}$ are obviously finite, and $T_{\ell}(B)$ is finitely
generated by Nakayama's lemma. Therefore, $T_{\ell}(B)$ is isomorphic to
$\Z_{\ell}^r$ for some nonnegative integer $r\le
\dim_{\F_{\ell}}(B_{\ell})$. Moreover, $T_{\ell}(B)=\{0\}$ if
and only if $B(\ell)$ is finite.
%because the projective limit of finite non-empty sets is also non-empty.

 For a field $k$ with separable closure $\ov k$ we
denote by $\Ga$ the Galois group $\Gal(\ov k/k)$. Let $X$ be a
geometrically integral smooth projective variety over $k$, and
let $\ov X=X\times_k\ov k$.

Let $\ell\ne \fchar(k)$ be a prime. Following  \cite[II, Sect.
3]{Gr} we recall that the exact Kummer sequence of sheaves in \'etale
topology
$$1\to\mu_{\ell^n}\to\G_m\to\G_m\to 1$$
gives rise to the (cohomological) exact sequence of Galois
modules
$$0\to \Pic (\ov X)/{\ell^n}\to\H_\et^2(\ov
X,\mu_{\ell^n}) \to \Br (\ov X)_{\ell^n}\to 0.$$
Since
$\Pic(\ov X)$ is an extension of the N\'eron--Severi group
$\NS (\ov X)$ by
a divisible group, we have $\Pic (\ov X)/{\ell^n}=\NS
(\ov X)/{\ell^n}$. We thus obtain the exact sequence of Galois
modules
\begin{equation}
0\to \NS(\ov X)/{\ell^n}\to\H_\et^2(\ov X,\mu_{\ell^n})
\to \Br (\ov X)_{\ell^n}\to 0. \label{e1}
\end{equation}
Since the groups $\H_\et^2(\ov X,\mu_{\ell^n})$ are finite,
the groups $\Br (\ov X)_{\ell^n}$ are finite as well \cite[II, Cor.
3.4]{Gr}.  On passing to the projective limit we get an exact
sequence of $\Ga$-modules
\begin{equation}
0\to \NS (\ov X)\otimes\Z_\ell\to\H_\et^2(\ov X,\Z_\ell(1))\to
T_\ell(\Br(\ov X))\to 0. \label{NS}
\end{equation}
Since  $T_\ell(\Br(\ov X))$ is a free $\Z_{\ell}$-module, this
sequence shows that the torsion subgroup of $\H_\et^2(\ov
X,\Z_\ell(1))$ is contained in $\NS (\ov X)\otimes\Z_\ell$, that is,
the torsion subgroups of $\H_\et^2(\ov X,\Z_\ell(1))$ and $\NS
(\ov X)\otimes\Z_\ell$ coincide, and so are both equal to 
$\NS (\ov X)(\ell)$.
%On the other hand, since $\NS
%(\ov X)\otimes\Z_\ell$ is finitely generated, the inclusion $\NS
%(\ov X)=\NS (\ov X)\otimes 1 \subset\NS (\ov X)\otimes\Z_\ell$
%induces a canonical isomorphism
%$$\NS (\ov X)(\ell) \cong (\NS (\ov X)\otimes\Z_{\ell})_\tors.$$
Tensoring the sequence (\ref{NS}) with $\Q_{\ell}$ (over
$\Z_{\ell}$), we get the exact sequence of $\Ga$-modules
\begin{equation}
0\to \NS (\ov X)\otimes\Q_\ell\to\H_\et^2(\ov X,\Q_\ell(1))\to
V_\ell(\Br(\ov X))\to 0, \label{Ql}
\end{equation}
where $V_\ell(\Br(\ov X))=T_\ell(\Br(\ov X))\otimes_{\Z_{\ell}}\Q_\ell$.
The Tate
conjecture for divisors \cite{tateW,tate,T91} asserts that if $k$
is finitely generated over its prime subfield, then
\begin{equation}
\H_\et^2(\ov X,\Q_\ell(1))^{\Ga}=\NS (\ov
X)^{\Ga}\otimes\Q_\ell.\label{TCD}
\end{equation}
\end{sect}
Note also that (\ref{e1}) gives rise to the exact sequence
of abelian groups
\begin{equation}
\begin{array}{c}
0\to (\NS (\ov X)/{\ell^n})^\Ga\to \H_\et^2(\ov
X,\mu_{\ell^n})^\Ga\to
\Br (\ov X)^\Ga_{\ell^n}\\
\\
\to \H^1(k,\NS (\ov X)/{\ell^n})\to
\H^1(k,\H_\et^2(\ov X,\mu_{\ell^n})).
\end{array} \label{e1.1}
\end{equation}

The lemma that follows is probably well known, cf.
\cite[Sect. 5, pp. 16--17 ]{KleimanM} and \cite[pp. 198--199 ]{DM}.

\begin{lem}
\label{hdg}
Let $L\in \NS(\ov X)^\Ga$ be a Galois invariant
hyperplane section class. Assume that $d=\dim(X)\ge 2$. If $\fchar(k)=0$,
then the kernel of the symmetric intersection pairing
$$\psi_0:\NS(\ov X) \times \NS(\ov X) \to \Z, \quad
x,y\mapsto x\cdot y\cdot L^{d-2},$$
is $\NS(\ov X)_{\tors}$. 

In any characteristic the same conclusion holds
under the following condition:

there exist a finite extension $k'/k$ with $k'\subset \ov k$,
and a prime $q\ne \fchar(k)$ such that $\Gal(\ov k/k')$ acts trivially on
$\NS (\ov X)$, the $\Gal(\ov k/k')$-module
$\H_\et^2(\ov X,\Q_q(1))$ is semisimple, and 
$\H_\et^2(\ov X,\Q_q(1))^{\Gal(\ov k/k')}=\NS (\ov X)\otimes\Q_q$.
\end{lem}
\begin{proof}
We start with the case of characteristic zero.
If $K$ is an algebraically closed field containing $k$, then
the N\'eron--Severi group $\NS(X\otimes_kK)$ is identified with
the group of connected components of the Picard scheme of $X$
\cite[Cor. 4.18.3, Prop. 5.3, Prop. 5.10]{Kleiman},
and so does not depend on $K$. Let $k_0 \subset k$
be a subfield finitely generated over $\Q$, over which
$X$ and $L$ are defined. Then
there exists a smooth projective variety
$X_0$ over the algebraic closure $\ov{k_0}$ of $k_0$ in $\ov k$,
such that $\ov X=X_0\times_{\ov{k_0}} \ov k$.
The natural map $\NS(X_0)\to \NS(\ov X)$ is bijective,
and therefore a group isomorphism.

For generalities on twisted classical cohomology groups
we refer the reader to
see \cite[Sect. 1]{DH} or \cite[Sect. 2.1]{DWI}. 

Fix an embedding $\ov{k_0}\hookrightarrow \C$ and consider the
complex variety $X_{\C}=X_0\times_{\ov{k_0}}\C$.
The natural map $\NS( X_0)\to \NS(X_{\C})$
is an isomorphism. Since the
intersection indices do not depend on the choice of an
algebraically closed ground field, it suffices to check the
non-degeneracy of $\psi_0$ for the complex variety $X_{\C}$.
In order to do so, consider the canonical embedding
$$\NS(X_{\C})\otimes\Q\hookrightarrow \H^2(X_{\C} (\C),\Q(1)),$$
and the symmetric bilinear form
$$
\rho : \ \H^2(X_{\C}(\C),\Q(1)) \times\ \H^2(X_{\C}(\C),\Q(1)) \to
\ \Q, \ \ \ x,y\mapsto x\cup y\cup L^{d-2} \label{p} .$$  The Hard
Lefschetz theorem says that the map
$$\H^2(X_{\C}(\C),\Q(1))\lra\H^{2d-2}(X_{\C}(\C),\Q(d-1)), 
\quad x\mapsto x\cup L^{d-2},$$
is an isomorphism of vector spaces over $\Q$. Poincar\'e
duality now implies that $\rho$ is non-degenerate. Let us show
that the restriction of $\rho$ to $\NS(\ov X)\otimes\Q \subset
\H^2(X_{\C}(\C),\Q(1))$ is also non-degenerate. Indeed, let
$P\subset \H^2(X_{\C}(\C),\Q(1))$ be the 
kernel of the multiplication by $L^{d-1}$. The
group $\H^2(X_{\C}(\C),\Q(1))$ is the orthogonal
direct sum $\Q L \oplus P$. On the one hand,
the form $\rho$ is positive definite
on $\Q L$ since $L$ is ample. On the other hand, the
restriction of $\rho$ to $P$ is negative definite, due to the
Hodge--Riemann bilinear relations \cite[Ch. V, Sect. 5, Thm.
5.3]{We}. This implies the non-degeneracy of $\rho$ on
$\NS(X_{\C})\otimes\Q$, because this space is the direct
sum of $\Q L$ and $(\NS(X_{\C})\otimes\Q) \cap P$.
To finish the proof, we note that the form induced
by $\rho$ on the N\'eron--Severi group coincides with $\psi_0$, whereas
the kernel of $\NS(X_{\C})\to \NS(X_{\C})\otimes\Q$ is the torsion
subgroup $\NS(X_{\C})_{\tors}=\NS(\ov X)_{\tors}$.
\medskip

Now let us prove the lemma in the case of arbitrary characteristic, assuming
the condition on the Galois module $\H_\et^2(\ov X,\Q_q(1))$. 

Let us replace $k$ by $k'$.
Consider the symmetric Galois-invariant $\Q_{q}$-bilinear form
$$
\rho_q : \ \H^2_\et(\ov{X},\Q_{q}(1))\  \times \
\H^2_\et(\ov{X},\Q_{q}(1))  \to \ \Q_{q}, \ \ \  x,y\mapsto x\cup y\cup
L^{d-2} \label{pl} .$$  The Hard Lefschetz theorem proved by
Deligne \cite{DW2} in all characteristics, says
that the map
$$h_L: \H^2_\et(\ov{X},\Q_{q}(1))\lra\H^{2d-2}_\et(\ov{X},\Q_{q}(d-1)),
\quad x\mapsto x\cup L^{d-2},$$
is an isomorphism of vector spaces over $\Q_{q}$. Thus $h_L$
is an isomorphism of Galois modules. Poincar\'e duality now
implies that $\rho_{q}$ is non-degenerate.

Since $h_L$ is an isomorphism of Galois modules, we have
 $$\H^{2d-2}_\et(\ov{X},\Q_{q}(d-1))^{\Ga}=h_L(\H^2_\et(\ov{X},\Q_{q}(1))^{\Ga})=
(\NS(\ov{X})\otimes\Q_{q})\cup L^{d-2}.$$ By the semisimplicity
of $\H^2_\et(\ov{X},\Q_{q}(1))$ there is a unique 
$\Ga$-invariant vector subspace $W$ which is also a
semisimple $\Ga$-submodule, such that
$$\H^2_\et(\ov{X},\Q_{q}(1))=(\NS(\ov{X})\otimes\Q_{q})\oplus W.$$
Our condition implies that $W^{\Ga}=\{0\}$. 
If $M\subset \H^2_\et(\ov{X},\Q_{q}(1))$ is a vector subspace which
is also a simple $\Ga$-submodule, and $M\to \Q_q$
is a non-zero $\Ga$-invariant linear form,
then $M$ is the trivial $\Ga$-module $\Q_q$.
It follows that the trivial $\Ga$-module
$\NS(\ov{X})\otimes\Q_{q}$ is orthogonal to $W$ with respect to 
$\rho_{q}$.
Now the non-degeneracy of $\rho_q$ implies that its restriction
$$\psi_{q}:\NS(\ov X)\otimes\Q_{q}\ \times\ \NS(\ov X)\otimes\Q_{q} \to \Q_{q},
 \quad x,y\mapsto x\cdot y\cdot L^{d-2},$$
is also non-degenerate. 
By the compatibility of the
cohomology class of the intersection of algebraic cycles and the
cup-product of their cohomology classes \cite[Ch. VI, Prop. 9.5
and Sect. 10]{MilnE}, the bilinear form $\psi_{q}$ is
obtained from $\psi_0$ by tensoring it with $\Q_{q}$. 
To finish the proof we note that
the kernel of $\NS(\ov X)\to \NS(\ov X)\otimes\Q_{q}$ is $\NS(\ov X)_{\tors}$.
\end{proof}

\begin{rema}
(i)
Since $\NS(\ov{X})$ is a finitely generated abelian group, there exists a finite extension $k'/k$ with $k'\subset \ov{k}$,
such that $\Gal(\ov{k}/k')$ acts trivially on $\NS(\ov{X})$. 

(ii)
Recall that $V:=\H_{\et}^2(\ov X,\Q_{\ell}(1))$ is a 
finite-dimensional $\Q_{\ell}$-vector space. 
Let $G_{\ell,k}$ be the image of 
$\Gamma=\Gal(\ov k/k)$ in $\Aut_{\Q_{\ell}}(V)$; 
it is a compact subgroup of $\Aut_{\Q_{\ell}}(V)$ and, 
by the $\ell$-adic version of Cartan's theorem \cite{SerreLie}, 
is an $\ell$-adic Lie subgroup of $\Aut_{\Q_{\ell}}(V)$.  
If $k'/k$ is a finite extension  with $k'\subset \ov{k}$,  
then $\Gamma'=\Gal(\ov k/k')$ is an open subgroup of finite index 
in $\Gamma$, hence the image $G_{\ell,k'}$ of $\Gamma'$ is 
an open subgroup of finite index in $G_{\ell,k}$. 
In particular, $G_{\ell,k}$ 
and $G_{\ell,k'}$ have the same Lie algebra, which is a $\Q_{\ell}$-Lie 
subalgebra of $\End_{\Q_{\ell}}(V)$. Applying Prop. 1 of 
\cite{SerreLocal}, we conclude that $V$  
is semisimple as a $G_{\ell,k'}$-module if and only if it is semisimple
as a $G_{\ell,k}$-module. It follows that 
$\H_{\et}^2(\ov{X},\Q_{\ell}(1))$ is semisimple 
as a $\Gamma'$-module if and only if it is semisimple as a $\Gamma$-module.
\end{rema}

The following statement was inspired by \cite[III, Sect. 8, pp.
143--147]{Gr} and \cite[Sect. 5]{tateB}.

\bpr \label{1} Let $X$ be a smooth projective geometrically
integral variety over a field $k$. Assume that one of the
following conditions holds.

\begin{itemize}
\item[(a)] $X$ is a curve or a surface.

\item[(b)] $\fchar(k)=0$.

\item[(c)] there exist a finite extension $k'/k$ with 
$k'\subset \ov k$,
and a prime $q\ne \fchar(k)$ such that $\Gal(\ov k/k')$ acts trivially on
$\NS (\ov X)$, the $\Gal(\ov k/k')$-module
$\H_\et^2(\ov X,\Q_q(1))$ is semisimple, and 
$\H_\et^2(\ov X,\Q_q(1))^{\Gal(\ov k/k')}=\NS (\ov X)\otimes\Q_q$.

\end{itemize}
Then for almost all primes $\ell$ the $\Ga$-module $\NS(\ov X)\otimes\Z_\ell$
is a direct summand of the $\Ga$-module $\H_\et^2(\ov X,\Z_\ell(1))$.
If {\rm (c)} is satisfied, then $\Br(\ov{X})^{\Ga}(q)$ is finite.
\epr

\begin{proof} (a)
If $X$ is a curve, then
$\H_\et^2(\ov X,\Z_\ell(1))=\NS(\ov X)\otimes\Z_\ell\cong\Z_\ell$, and 
there is nothing to prove. Note that in this case $\Br(\ov{X})=0$
\cite[III, Cor. 5.8]{Gr}. Thus from now on we assume that $\dim(X)\ge 2$.

Let $X$ be a 
surface, $n=|\NS(\ov X)_{\tors}|$. The cycle map defines the
commutative diagram of pairings given by the intersection index
and the cup-product:
\begin{equation}
\begin{array}{ccccc}
\H^2_\et(\ov X,\Z_{\ell}(1))&\times&
\H^2_\et(\ov X,\Z_{\ell}(1))&\to&\Z_\ell\\
\uparrow&&\uparrow&&\uparrow\\
\NS(\ov X)&\times&\NS(\ov X)&\to&\Z
\end{array} \label{p3}
\end{equation}
The diagram commutes by the compatibility of the
cohomology class of the intersection of algebraic cycles and the
cup-product of their cohomology classes \cite[Ch. VI, Prop. 9.5
and Sect. 10]{MilnE}.
 The kernel of the pairing on the N\'eron--Severi group is its torsion
 subgroup.
%moreover, the upper pairing is a perfect auto-duality of the
%Galois module $\H^2_\et(\ov X,\Z_{\ell}(1))/_\tors$.
Let $\delta$ be the discriminant of the induced bilinear form on 
$\NS(\ov X)/_\tors$. Let
$H$ be the $\Ga$-module $H^2_\et(\ov X,\Z_{\ell}(1))/_\tors$,
and let $\psi$ be the Galois-invariant $\Z_{\ell}$-bilinear form
on $H$ coming from the top pairing of (\ref{p3}). Let
$N$ be the $\Ga$-submodule $\NS(\ov X)/_\tors \otimes \Z_{\ell} \subset H$.
It is clear that $N$ is
a free $\Z_{\ell}$-submodule of $H$,
and $\delta$ is the discriminant of the restriction of
$\psi$ to $N$. Let $N^\perp$ be the orthogonal
complement to $N$ in $H$ with respect to $\psi$; $N^\perp$ is obviously a
$\Ga$-submodule of $H$.

Applying Lemma
\ref{discr} (with $\Lambda=\Z_{\ell}$) we conclude that
$$N\cap N^\perp=\{0\}\quad  \text{and}\quad  \delta^2 H\subset N\oplus N^\perp.$$
Now let $\tilde{M}$ be the preimage of $N^\perp$ in $\H^2_\et(\ov
X,\Z_{\ell}(1))$. Clearly, $\tilde{M}$ is a Galois submodule, and
$\tilde{M}\cap (\NS(\ov X)\otimes\Z_\ell)$ is the torsion subgroup
of $\NS(\ov X)\otimes\Z_\ell$, and therefore coincides with
$\NS(\ov X)(\ell)$. It is also clear that
$$\delta^2\, \H^2_\et(\ov X,\Z_{\ell}(1))\subset (\NS(\ov X)\otimes\Z_\ell)+
\tilde{M}.$$ Let us put $M=n\tilde M\subset\tilde M$. We have
$$M\cap (\NS(\ov X)\otimes\Z_\ell)=\{0\} \quad  \text{and}\quad
n\delta^2\, \H^2_\et(\ov X,\Z_{\ell}(1))\subset (\NS(\ov
X)\otimes\Z_\ell)\oplus M.$$
Since $\H^2_\et(\ov X,\Z_{\ell}(1))$
is a finitely generated $\Z_{\ell}$-module, 
$(\NS(\ov X)\otimes\Z_\ell)\oplus M$ is a subgroup of finite index 
in $\H^2_\et(\ov X,\Z_{\ell}(1))$. 
This index is 1 if $\ell$ does not divide $n\delta$.

(b) and (c). Let us choose a $\Ga$-invariant hyperplane section class
$L\in \NS(\ov X)^\Ga$. By Lemma \ref{hdg} the
symmetric bilinear form on $\NS(\ov X)/_{\tors}$
induced by $\psi_0$ is non-degenerate.
Let $\delta\in \Z$ be the discriminant of this form, $\delta\not=0$.
Let us consider the Galois-invariant
symmetric $\Z_{\ell}$-bilinear form
$$\psi_1:\H^2_\et(\ov X,\Z_{\ell}(1))\times \H^2_\et(\ov
X,\Z_{\ell}(1))\to\Z_\ell, \quad x,y\mapsto x\cup y\cup L^{d-2}.$$ 
The compatibility of (the cohomology class of) the intersection of
algebraic cycles and the cup-product of their cohomology classes
 \cite[Ch. VI, Prop. 9.5 and Sect. 10]{MilnE} implies that
 the restriction of $\psi_1$ to
$\NS(\ov X)\otimes\Z_\ell$ coincides with the form induced by $\psi_0$.
It follows from the Hard Lefschetz theorem and Poincar\'e duality
that $\ker(\psi_1)=\H^2_\et(\ov X,\Z_{\ell}(1))_{\tors}$.

Let $H$ be the $\Ga$-module $\H^2_\et(\ov X,\Z_{\ell}(1))/_\tors$,
and let $\psi$ be the Galois-invariant $\Z_{\ell}$-bilinear form on $H$
defined by $\psi_1$. Let
$N$ be the $\Ga$-submodule $\NS(\ov X)/_\tors \otimes \Z_{\ell} \subset H$.
It is clear that $N$ is a free $\Z_{\ell}$-submodule of $H$,
and the discriminant of the restriction of
$\psi$ to $N$ is $\delta$. 
The rest of the proof is the same as in case (a).

Now suppose that under the condition of (c) the group
 $\Br(\ov{X})^{\Ga}(q)$ is infinite. 
Since $\Br(\ov{X})^{\Ga}\subset \Br(\ov{X})^{{\rm Gal}(\ov k/k')}$,
we can extend the ground field from $k$ to $k'$.
 For any $n$ the group
 $\Br(\ov{X})_{q^n}$ is finite, thus there is an element
 of order $q^n$ in $\Br(\ov{X})^{\Ga}_{q^n}$
for every $n$, i.e., the set $S(n)$ of elements of
 order $q^n$ in  $\Br(\ov{X})^{\Ga}_{q^n}$ is non-empty for all $n$.
 Since the projective
 limit of non-empty finite sets $S(n)$ is a non-empty subset of
 $T_{q}(\Br(\ov{X})^{\Ga})\setminus \{0\}$, we conclude
 that
$$T_{q}(\Br(\ov{X}))^{\Ga}=T_{q}(\Br(\ov{X})^{\Ga})\ne \{0\}.$$
It follows that $V_q(\Br(\ov X))^{\Ga}\ne \{0\}$.
However,
the semisimplicity of $\H^2_\et(\ov{X},\Q_{q}(1))$ 
implies that the exact sequence
of Galois modules (\ref{Ql}) splits, that is, 
$$\H^2_\et(\ov{X},\Q_{q}(1))\cong
(\NS(\ov{X})\otimes \Q_{q})\oplus V_q(\Br(\ov X)).$$ 
By condition (c) we have
$V_q(\Br(\ov X))^{\Ga}= \{0\}$. This contradiction
proves the finiteness of $\Br(\ov{X})^{\Ga}(q)$.
\end{proof}

\bco \label{c} Let $X$ be a smooth projective geometrically
integral variety over a field $k$. Assume that $X/k$
satisfies one
of the conditions {\rm (a)}, {\rm (b)}, {\rm (c)} of Proposition $\ref{1}$.
Then the map
$\H^1(k,\NS (\ov X)\otimes\Z/{\ell})\to \H^1(k,\H_\et^2(\ov
X,\mu_{\ell}))$ in $(\ref{e1.1})$ is injective for almost all
$\ell$. \eco

\begin{proof}
By  Proposition \ref{1}, the
$\Ga$-module $\NS (\ov X)/{\ell}=(\NS (\ov X)\otimes\Z_{\ell})/\ell$
is a direct summand of the
$\Ga$-module $\H^2_\et(\ov X,\Z_\ell(1))/\ell$ for almost all $\ell$. 
 We have an exact sequence
$$0\to \H^2_\et(\ov X,\Z_\ell(1))/\ell\to
\H^2_\et(\ov X,\mu_\ell)\to \H^3_\et(\ov X,\Z_\ell(1))_\ell\to 0.$$ 
By a theorem of Gabber \cite{G}, for almost all
$\ell$ the $\Z_{\ell}$-module $\H^3_\et(\ov X,\Z_\ell)$ has
no torsion. Since $\H^3_\et(\ov
X,\Z_\ell)$ and $\H^3_\et(\ov X,\Z_\ell(1))$ are isomorphic
as abelian groups, for
almost all $\ell$ we have $\H^3_\et(\ov
X,\Z_\ell(1))_\ell=\{0\}$, hence $\H^2_\et(\ov
X,\mu_\ell)=\H^2_\et(\ov X,\Z_\ell(1))/\ell$. Thus
$\NS (\ov X)/{\ell}$ 
is a direct summand of $\H^2_\et(\ov X,\mu_\ell)$. 
This proves the corollary.
\end{proof}

\bco \label{2.3} Suppose that $k$ is finitely generated over its
prime subfield, and $\fchar(k)\ne 2$. Let $A$ be an abelian variety
over $k$. Then for all $\ell\not=\fchar(k)$ the subgroup 
$\Br (\ov A)^\Ga(\ell)$ is finite. 
\eco

\begin{proof}
Let $\ell$ be a prime different from $\fchar(k)$.
The Tate conjecture for divisors (\ref{TCD})
is true for any abelian variety $A$ over such a field; in addition, 
the natural Galois action on the
$\ell$-adic cohomology groups of $\ov A$ is semisimple. (These
assertions were proved by the second named author \cite{Z75,Z76}
in finite characteristic not equal to $2$, and by Faltings
\cite{F,F86} in characteristic zero.) This implies that $A$
satisfies condition (c) of Proposition \ref{1} for every prime
$q\ne \fchar(k)$. Now the result follows from the last assertion
of Proposition \ref{1}.
\end{proof}

\section{Proof of Theorem \ref{t1}}
%{Applying the Tate conjecture with finite coefficients}
\label{TateC}

Let $A$, $A'$ be abelian varieties over an {\sl arbitrary} field $k$.
We write $\Hom(A,A')$ for the group of homomorphisms $A\to A'$.
We have 
$$\Hom(A,A')=\Hom_\Ga(\ov A,\ov A')=\Hom(\ov A,\ov A')^\Ga.$$
Since $\Hom(\ov A,\ov A')$ has no torsion,
$\Hom(A,A')/n$ is a subgroup of $\Hom(\ov A,\ov A')/n$.

Let $A^t$ be the dual abelian variety of $A$. We have $(A^t)^t=A$
(\cite[Ch. V, Sect. 2, Prop. 9]{Lang}, \cite[p. 132]{Mu}). Every
divisor $D$ on $\ov A$ defines the homomorphism $\ov A\to \ov A^t$
sending $a\in A(\ov k)$ to the linear equivalence class of
$T_a^{*}(D)-D$ in $\Pic^0(\ov A)$, where $T_a$ is the translation
by $a$ in $A$. If $L$ is the algebraic equivalence class of $D$ in
$\NS(\ov A)$, then this map depends only on $L$, and is denoted by
$\phi_L:\ov A \to \ov A^t$ \cite[Sect. 8]{Mu}. For 
$\alpha\in \Hom(\ov A,\ov A^t)$ 
we denote by $\alpha^t\in \Hom(\ov A,\ov A^t)$ 
the transpose of $\alpha$. Note that $\phi_L^t=\phi_L$.
Moreover, if we set
$$\Hom(\ov A,\ov A^t)_\sym:=\{u\in \Hom(\ov A,\ov A^t)\mid
u=u^t\},$$ then the group homomorphism
$$\NS(\ov A)\to \Hom(\ov A,\ov A^t)_\sym, \quad L\mapsto \phi_L,$$
is an isomorphism \cite{Lang}, \cite[Sect. 20, formula (I) and Thm.
1 on p. 186, Thm. 2 on p. 188 and Remark on p. 189]{Mu}. For any
$\alpha\in \Hom(\ov A,\ov A^t)$ we have
$(\alpha^t)^t=\alpha$, and thus
\begin{equation}
\alpha+\alpha^t\in \Hom(\ov A,\ov A^t)_\sym. \label{symm}
\end{equation}

\begin{sect}
\label{skew} Let $\ell$ be a prime different from the
characteristic of $k$, and $i$ a positive integer,
$n=\ell^i$. The kernel $A_{n}$ of the multiplication by $n$ in
$A(\ov{k})$ is a Galois submodule, isomorphic to $(\Z/n)^{2\dim(A)}$
as an abelian group.

The natural map
$\Hom(\ov A,\ov A')/n\to \Hom(A_{n},A'_n)$ 
is {\sl injective} \cite[p. 124]{MilneAV}. It
commutes with the Galois action on both sides; in
particular, the image of $\Hom({A},A')/n\subset \Hom(\ov A,\ov A')/n$ 
lies in $\Hom_{\Gamma}(A_{n},A'_n)$.

For any $\alpha\in \Hom(\ov A,\ov A^t)$ and any $x,y\in A_n$ 
we have (\cite[Ch. VII, Sect. 2, Thm. 4]{Lang}, \cite[p. 186]{Mu})
$$e_{n,A^t}(\alpha x,y) =e_{n,A}(x,\alpha^t{y}).$$ 
Thus $\Hom(\ov A,\ov A^t)_\sym/n$ is a subgroup of 
$$\Hom (A_{n},A^t_{n})_\sym:=\{u \in \Hom
(A_{n},A^t_{n})\mid e_{n,A^t}(u x,{y})=e_{n,A}(x,uy), \
\forall x,y\in A_n\}.$$ 
Moreover, if $\ell$ is odd, then, by (\ref{symm}), we have
\begin{equation}
\Hom(\ov A,\ov A^t)_\sym/n=\Hom(\ov A,\ov A^t)/n\cap\Hom
(A_{n},A^t_{n})_\sym.
\label{homsym}
\end{equation}

\begin{rema}
\label{skewsym}
The two (non-degenerate, Galois-equivariant) Weil pairings
$$e_{n,A}:A_{n}\times A^t_{n}\to\mu_{n} \ \ \text{and} \ \
e_{n,A^t}:A^t_{n}\times A_{n}\to\mu_{n}$$ differ by $-1$ \cite[Ch.
VII, Sect. 2, Thm. 5(iii) on p. 193]{Lang}, that is, 
$$e_{n,A^t}(y,x)=-e_{n,A}(x,y)$$ 
for all $x\in A_n$, $y\in A^t_{n}$.
 Since for each $u\in \Hom(A_{n},A^t_{n})$ we have
$$e_{n,A}(x,uy)=-e_{n,A^t}(uy,{x})=-e_{n,A}(y,u^t{x}),$$
we conclude that $u$ lies in $\Hom (A_{n},A^t_{n})_\sym$ if and
only if the bilinear form $e_{n,A}(x,uy)$ is
{\sl skew-symmetric}, that is, for any $x,y\in A_n$ we have
$$
e_{n,A}(x,uy)=-e_{n,A}(y,ux).
$$
\end{rema}
\end{sect}
\begin{sect}
\label{tateL}
For a module $M$ over a commutative ring $\Lambda$ we denote by
$\mathrm{S}^2_{\Lambda} M$ the submodule of $M\otimes_{\Lambda}M$ 
generated by $x\otimes x$ for all $x\in M$. 
Let $\wedge^2_{\Lambda}M=
(M\otimes_{\Lambda}M)/\mathrm{S}^2_{\Lambda} M$. We have 
$x\otimes y+y\otimes x\in \mathrm{S}^2_{\Lambda} M$; these
elements generate $\mathrm{S}^2_{\Lambda} M$ if $2$ is invertible in $\Lambda$.

 {From} the Kummer sequence one obtains the well known
canonical isomorphism $\H^1_\et(\ov
A,\mu_{n})=\Pic(\ov{A})_{n}=A^t_{n}$. 
Thus we have
canonical isomorphisms of Galois modules (cf. \cite[Sect. 2]{Ber},
\cite{MilnE}, \cite{MilneAV}):
$$
\H^2_\et(\ov A,\mu_{n})=\wedge_{\Z/n}^2
A^t_{n}(-1)=\Hom(\wedge_{\Z/n}^2 A_{n},\mu_n). 
$$
Clearly, there is a canonical embedding of Galois modules
$$\Hom(\wedge_{\Z/n}^2 A_{n},\mu_n)\hookrightarrow
\Hom(A_{n},A^t_{n}),$$ whose image coincides with the set
of $u:A_{n}\to A^t_{n}$ such that the bilinear form $
e_{n,A}(x,uy)$ is {\sl alternating}, i.e., $e_{n,A}(x,ux)=0$ for all
$x\in A_n$. Combining it with Remark \ref{skewsym}, we conclude
that if $\ell$ is {\sl odd}, then there is a canonical isomorphism
of Galois modules
\begin{equation}
\H^2_\et(\ov A,\mu_{n})\cong \Hom (A_{n},A^t_{n})_\sym .
\label{H2sym}
\end{equation}

\end{sect}

 Let us recall a variant of the Tate conjecture on
homomorphisms that first appeared in \cite{Z77}. 

\bpr \label{TLB}
Let $k$ be a field finitely generated over its prime subfield,
$\fchar(k)\ne 2$. If $A$ and $A'$ are abelian varieties over $k$,
then the natural injection
\begin{equation}
\Hom(A,A')/\ell \hookrightarrow
\Hom_\Ga(A_{\ell},A'_{\ell}) \label{9}
\end{equation}
is an isomorphism for almost all $\ell$. 
\epr
\begin{proof}
In the finite characteristic case this is
proved in \cite[Thm. 1.1]{Z77}. When $A=A'$ and $k$ is a
number field, Cor. 5.4.5 of \cite{Z85} (based on the
results of Faltings \cite{F}) says that for almost all $\ell$
we have
\begin{equation}
\End(A)/\ell =\End_\Ga(A_{\ell}). \label{10}
\end{equation}
The same proof works over arbitrary fields that are finitely
generated over $\Q$, provided one replaces the reference to Prop.
3.1 of \cite{Z85} by the reference to the corollary on p. 211 of
Faltings \cite{F86}. Applying (\ref{10}) to the abelian variety
$A\times A'$, we deduce that (\ref{9}) is a bijection.
\end{proof}

\ble \label{2.2} Let $k$ be a field finitely generated over its
prime subfield, $\fchar(k)\ne 2$, and let $A$ be an abelian variety
over $k$.
Then for almost all $\ell$ we have the following statements:
\begin{itemize}
\item[(i)] the injective map $(\NS (\ov
A)/\ell)^\Ga\hookrightarrow \H_\et^2(\ov A,\mu_\ell)^\Ga$
in $(\ref{e1.1})$ is an isomorphism;

%\item[(ii)] the natural map $\NS(\ov A)^\Ga/\ell\to (\NS
%(\ov A)/\ell)^\Ga$ is an isomorphism;

\item[(ii)]  $\Br(\ov X)^\Ga(\ell)=\{0\}$.
\end{itemize}
\ele
\begin{proof}
Suppose that $\ell$ is odd.
By (\ref{homsym}) we have
$$\Hom(\ov A,\ov A^t)_\sym/\ell=\Hom(\ov A,\ov A^t)/\ell\cap\Hom
(A_{\ell},A^t_{\ell})_\sym.$$
Proposition \ref{TLB} implies that for almost all $\ell$ we have
$$
\Hom(A,A^t)/\ell=\Hom(A_\ell,A^t_\ell)^{\Ga}=\Hom_{\Ga}(A_\ell,A^t_\ell).
$$
We thus obtain an isomorphism
\begin{equation}
\Hom(A,A^t)_{\sym}/\ell= \Hom_\Ga(A_\ell,A^t_\ell)_{\sym}.\label{aaa}
\end{equation}
The left hand side of (\ref{aaa}) is
$\Hom_\Ga(\ov A,\ov A^t)_\sym/\ell\cong\NS(\ov A)^\Ga/\ell$,
see the beginning of this section. The
right hand side of (\ref{aaa})
is isomorphic to $\H_\et^2(\ov A,\mu_\ell)^\Ga$
by (\ref{H2sym}). It follows that $\NS(\ov A)^\Ga/\ell$ and 
$\H_\et^2(\ov A,\mu_\ell)^\Ga$ have the same number of elements. 
Since $\NS(\ov A)$ has no torsion,
$\NS(\ov A)^\Ga/\ell$ is a subgroup of $(\NS(\ov A)/\ell)^\Ga$,
and hence the injective map in (i) is bijective. Statement
(ii) follows from (i), Corollary \ref{c} and the exact
sequence (\ref{e1.1}).
\end{proof}

\medskip

\begin{proof}[End of proof of Theorem $\ref{t1}$] Let $A$ be an abelian
variety over $k$, and $X$ a principal homogeneous space of $A$. 
In characteristic 0 (resp. in characteristic $p$)
it suffices to show that $\Br(\ov X)^\Ga$
(resp. $\Br(\ov X)^\Ga ({\non}p)$) is finite. For this we can go
over to a finite extension $k'/k$ such that
$X\times_kk'\simeq A\times_kk'$, and 
so assume that $X=A$.
The theorem now follows from Lemma \ref{2.2} (ii) and Corollary
\ref{2.3}. \end{proof}

\section{Proof of Theorem \ref{t}}
\label{K33}

\begin{sect}
\label{H1A}
In this subsection we recall some well known results which  
will be used later in this section.

Let $A$ be an abelian variety over a field $k$, $\ell$
a prime different from $\fchar(k)$, $n=\ell^i$.
Let $\pi_1^{\et}(\ov{A},0)^{(\ell)}$ be the maximal abelian
$\ell$-quotient of the Grothendieck \'etale fundamental group
$\pi_1^{\et}(\ov{A},0)$.
 Let us consider the Tate $\ell$-module
$T_{\ell}(A):=T_{\ell}(A(\ov{k}))$.
It is well known \cite{Lang,Mu} that $T_{\ell}(A)$ is a
free $\Z_{\ell}$-module of rank $2\dim(A)$ equipped with a
natural structure of a $\Ga$-module, and the natural map
$T_{\ell}(A)/n \to A_n$
is an isomorphism of Galois modules. Recall \cite[pp.
129--130]{MilneAV} that the isogeny
$\ov{A}\stackrel{n}{\to}\ov{A}$ is a Galois \'etale covering with
the Galois group $A_n$ acting by translations. 
This defines a canonical surjection
$f_n:\pi_1^{\et}(\ov{A},0)^{(\ell)}\twoheadrightarrow A_n$. The $f_n$
glue together into a canonical isomorphism of
Galois modules
$\pi_1^{\et}(\ov{A},0)^{(\ell)}\to T_{\ell}(A)$, which 
induces  the canonical isomorphisms of Galois modules
$$\H^1_{\et}(\ov{A},\Z_{\ell})=\Hom_{\Z_{\ell}}(\pi_1^{\et}(\ov{A},0)^{(\ell)},\Z_{\ell})=
\Hom_{\Z_{\ell}}( T_{\ell}(A),\Z_{\ell}).$$ 
Since $\H^j_{\et}(\ov{A},\Z_{\ell})$ is torsion-free for any 
$j$ \cite[Thm. 15.1(b) on p. 129]{MilneAV}, the reduction
modulo $n$ gives rise to natural isomorphisms of Galois modules
$$\H^1_{\et}(\ov{A},\Z/n)=\H^1_{\et}(\ov{A},\Z_{\ell})/n=
\Hom(A_{\ell},\Z/n).\label{H1n}$$

Now suppose that we are given a field embedding
$\ov{k}\hookrightarrow\C$. Let us consider the complex abelian
variety $B=A(\C)$. The exponential map establishes a
canonical isomorphism of compact Lie groups
$\Lie(B)/\Pi \to B$
\cite[Sect. 1]{Mu}. Here $\Lie(B)\cong \C^{\dim(B)}$ is the
tangent space to $B$ at the origin, $\Pi$ is a discrete lattice of
rank $2\dim(B)$, and the natural map
$\H_1(B,\Z)\otimes\R\to\Lie(B)$
is an isomorphism of real vector spaces. Clearly, $V$ is the
universal covering space of $B$, and the fundamental group
$\pi_1(B,0)=\H_1(B,\Z)=\Pi$ acts on $V$ by translations.
We have
$$B_n=\frac{1}{n}\Pi/\Pi\subset V/\Pi=B.$$
%There is a natural comparison group homomorphism
%\begin{equation}
%\pi_1(B,0)\to
%\pi_1^{\et}(B,0)\twoheadrightarrow\pi_1^{\et}(B,0)^{(\ell)}.\label{fund}
%\end{equation}
The isogeny $B\stackrel{n}{\to}B$ is an unramified Galois covering of
connected spaces (in the classical topology) with the Galois group $B_n$,
corresponding to the subgroup $n\Pi\subset \Pi$. It is identified with
$V/n\Pi\to V/\Pi$, and the corresponding
homomorphism $\varphi_n:\Pi \twoheadrightarrow
B_n=\frac{1}{n}\Pi/\Pi$ sends $c$ to $\frac{1}{n}c+\Pi$. 
The comparison theorem for fundamental groups implies that
$\varphi_n$ coincides with the composition
$$\pi_1(B,0)\to \pi_1^{\et}(B,0)\to
\pi_1^{\et}(B,0)^{(\ell)}\stackrel{f_n}{\longrightarrow} B_n.$$
We obtain the following sequence of homomorphisms
\begin{equation}
\begin{array}{c}
\Hom(B_n,\Z/n)\hookrightarrow\Hom(\pi_1^{\et}(B,0)^{(\ell)},\Z/n)
=\\
\\
 \Hom(\pi_1^{\et}(B,0),\Z/n)\to \Hom(\pi_1(B,0),\Z/n).
\end{array}
\label{fundn}
\end{equation}
The same comparison theorem implies
that the last map in (\ref{fundn}) is bijective. It follows easily
that all the homomorphisms in (\ref{fundn}) are isomorphisms.
Recall that
$$\Hom(\pi_1^{\et}(B,0)^{(\ell)},\Z/n)=\H^1_{\et}(B,\Z/n),
\quad \Hom(\pi_1(B,0),\Z/n)=\H^1(B,\Z/n).$$ Note also that $\varphi_n$
establishes a canonical isomorphism
$$\Pi/n=\pi_1(B,0)/n \to B_n, \quad c\mapsto
\frac{1}{n}c+\Pi,$$ which gives us the canonical isomorphisms
$$\H^1(B,\Z/n)=\H^1(B,\Z)/n=\Hom(B_n,\Z/n)=\H^1_{\et}(B,\Z/n).$$
Taking the projective limits with respect to $i$ (recall that
$n=\ell^i$), we get  the canonical isomorphisms
$$\H^1(B,\Z)\otimes \Z_{\ell}=\Hom_{\Z_\ell}(T_{\ell}(B),\Z_{\ell})
=\H^1_{\et}(B,\Z_{\ell}).$$
On the other hand, taking the projective limit of the $\varphi_n$,
we get  the natural map \cite[p. 237]{Mu}
$$\H_1(B,\Z)=\Pi \to T_{\ell}(B), \quad x\mapsto
\{x/\ell^i\}_{i=1}^{\infty},$$ which extends by
$\Z_{\ell}$-linearity to the natural isomorphism of
$\Z_{\ell}$-modules
$$\varphi^{(\ell)}:\H_1(B,\Z)\otimes \Z_{\ell}=\Pi\otimes\Z_{\ell}\cong T_{\ell}(B).$$
We have
\begin{equation}
A_n=B_n=\H_1(B,\Z)/n.\label{K3n}
\end{equation}
The comparison theorem for \'etale and classical cohomology implies 
$\H^1_{\et}(\ov{A},\Z/n)=\H^1(B,\Z/n)$, thus
we obtain
\begin{equation}
\begin{array}{c}
\H^1_{\et}(\ov{A},\Z/n)=\Hom(A_n,\Z/n)=\Hom(B_n,\Z/n)=
\Hom(\H_1(B,\Z)/n,\Z/n), \\
\\
T_{\ell}(A)=T_{\ell}(B)=\H_1(B,\Z)\otimes\Z_{\ell},\\
\\
\Hom_{\Z_{\ell}}(T_{\ell}(A),\Z_{\ell})=
\H^1_{\et}(\ov{A},\Z_{\ell})=
\H^1_{\et}(B,\Z_{\ell})=\Hom_{\Z_{\ell}}(T_{\ell}(B),\Z_{\ell}).
\end{array}
\label{K3compT}
\end{equation}
\end{sect}

\bigskip

\ble \label{3.1} Let $M$ and $N$ be subgroups of $\Z^n$ such that
$M\cap N=0$.
Then for almost all $\ell$ the natural maps
$M/\ell\to (\Z/\ell)^n$ and $N/\ell\to
(\Z/\ell)^n$ are injective, and the intersection of their images is
$\{0\}$.
\ele

\begin{proof}
There is a subgroup $L\subset\Z^n$ such that $L\cap
(M\oplus N)=0$, and $L\oplus M\oplus N$ is of finite index in
$\Z^n$. For all $\ell$ not dividing this index, the canonical map
$M/\ell\to (\Z/\ell)^n$ and the similar map for $N$ are
injective. Moreover, $(\Z/\ell)^n$ is the direct sum of
$L/\ell$, $M/\ell$ and $N/\ell$. This
proves the lemma. \end{proof}
\medskip

\ble \label{K3} Let $X$ be a K3 surface over a field $k$ finitely
generated over $\Q$. 
Then the injective map 
$(\NS (\ov X)/\ell)^\Ga\to\H^2_\et(\ov X,\mu_\ell)^\Ga$ in
$(\ref{e1.1})$ is an isomorphism for almost all primes $\ell$. \ele

\begin{proof} It suffices to prove the lemma for 
a finite extension $k'/k$, $k'\subset \ov k$, and
$\Ga'=\Gal(\ov k/k')\subset \Ga$.
Indeed, for any $\Ga$-module $M$ the composition
of the natural inclusion $M^\Ga\hookrightarrow M^{\Ga'}$ and the
norm map $M^{\Ga'}\to M^\Ga$ is the multiplication by the degree $[k':k]$.
Hence if $(\NS (\ov X)/\ell)^{\Ga'}\to\H^2_\et(\ov X,\mu_\ell)^{\Ga'}$
is surjective for all primes $\ell$ not dividing a certain
integer $N$, then 
so is the original map $(\NS (\ov X)/\ell)^\Ga\to\H^2_\et(\ov X,\mu_\ell)^\Ga$
for all primes $\ell$ not dividing $N[k':k]$.
In particular, we can assume without loss of generality
that $\Ga$ acts trivially on $\NS (\ov X)$.

Now let us fix an embedding $\ov{k}\hookrightarrow\C$ and
identify $\ov{k}$ with its image in $\C$.

The group
$\H^2(X(\C),\Z(1))\simeq\Z^{22}$ has a natural $\Z$-valued
bilinear form $\psi$ given by the intersection index. By 
Poincar\'e duality $\psi$ is {\sl unimodular}, i.e., the map
$\H^2(X(\C),\Z(1))\to \Hom(\H^2(X(\C),\Z(1)),\Z)$ induced by
$\psi$ is an isomorphism. Since $X(\C)$ is  simply-connected we have
$\H^1(X(\C),\Z)=\{0\}$, and by Poincar\'e duality this implies
$\H^3(X(\C),\Z)= \{0\}$. Recall that
$\NS (\ov X)=\NS(X_{\C})$ (see the beginning of the proof of 
Lemma \ref{hdg}). Since $X(\C)$ is simply-connected we have
$$\Pic(X_{\C})=\NS(X_{\C})=\NS (\ov X)=\Pic(\ov X).$$

We define the lattice of transcendental cycles $T_X$ as the
orthogonal complement to the injective image of $\NS (\ov X)$ in
$\H^2(X(\C),\Z(1))$. The restriction of $\psi$ to $\NS(\ov X)$ is
non-degenerate, and we write $\delta$ for the absolute value of
the corresponding discriminant. Then $\NS(\ov X)\cap T_X=0$, and
$\NS(\ov X)\oplus T_X$ is a subgroup of $\H^2(X(\C),\Z(1))$ of
finite index $\delta$. 
Let $\ell$ be a prime not dividing $\delta$. Then we have
$$\H^2(X(\C),\Z(1))/\ell=(\NS(\ov X)/\ell)\oplus
(T_X/\ell).$$ The restriction of the $\Z/\ell$-valued
pairing induced by $\psi$ to $\NS(\ov X)/\ell$ is a
non-degenerate $\Z/\ell$-bilinear form, so that
$T_X/\ell$ is the orthogonal complement to $\NS(\ov
X)/\ell$. Since $\H^3(X(\C),\Z)=\{0\}$, we have
$\H^2(X(\C),\Z(1))/\ell=\H^2(X(\C),\mu_{\ell})$. The
comparison theorem gives an isomorphism of $\Z_{\ell}$-modules
\begin{equation}
\H^2_\et(\ov
X,\Z_\ell(1))=\H^2(X(\C),\Z(1))\otimes\Z_{\ell},\label{Zl}
\end{equation}
which is compatible with cup-products \cite[Prop. 6.1, p.
197]{DM}, \cite[Example 2.1(b), pp. 28--29]{DH}. Reducing
modulo $\ell$ we get an isomorphism of $\Z/\ell$-vector
spaces $\H^2_\et(\ov X,\mu_\ell)=\H^2(X(\C),\Z(1))/\ell$,
compatible with cup-products. Thus for $\ell$ not dividing
$\delta$ we have an orthogonal direct sum
$$\H^2_\et(\ov X,\mu_\ell)=(\NS (\ov X)/\ell)\oplus
(T_X/\ell),$$ so that for these $\ell$ the abelian
group $T_X/\ell$ carries a natural $\Ga$-(sub)module
structure. (Here we use the compatibility of the cycle maps
$\Pic(\ov X)\to  \H^2_\et(\ov X,\mu_\ell)$ and 
$\Pic(\ov X)\to \H^2(X(\C),\mu_{\ell})$, see
 \cite[Prop. 3.8.5, pp. 296--297]{SGA5}.)

Let $L\in\Pic (\ov X)=\NS (\ov X)$ be a $\Ga$-invariant
hyperplane section class, and $P\subset\H^2(X(\C),\Z(1))$ the
orthogonal complement to $L$ with respect to $\psi$. Then
(\ref{Zl}) implies that $P_{\ell}=P\otimes\Z_{\ell}$ is both a
Galois and $\Z_{\ell}$-submodule of $\H^2_\et(\ov
X,\Z_\ell(1))$. It is clear that $P_{\ell}$ is
the orthogonal
complement to $L$ in $\H^2_\et(\ov X,\Z_\ell(1))$ with respect to
the Galois-invariant intersection pairing
$$\psi_{\ell}:\H^2_\et(\ov X,\Z_\ell(1))\times \H^2_\et(\ov
X,\Z_\ell(1))\to\Z_{\ell}.$$ Similarly,
$T_X\otimes\Z_{\ell}$ is the orthogonal complement to
$\NS(\ov{X})\otimes\Z_{\ell}$ in $\H^2_\et(\ov X,\Z_\ell(1))$
with respect to $\psi_{\ell}$, and so is a Galois submodule.

Let $C^+(P)$ be the even Clifford $\Z$-algebra of $(P,\psi)$.
The complex vector space $P_\C:=P\otimes\C$ inherits
from $\H^2(X(\C),\C(1))$ the Hodge
decomposition of type $\{(-1,1),(0,0),(1,-1)\}$ with Hodge numbers
$h^{1,-1}=h^{-1,1}=1$. By the
Lefschetz theorem, $T_X$ intersects trivially with the
$(0,0)$-subspace. The $\Z$-algebra $C^+(P)$ naturally carries a
weight zero Hodge structure of type
$\{(-1,1),(0,0),(1,-1)\}$ induced by the Hodge structure on $P$
(via the identification $C^+(P)=\oplus_i \wedge^{2i}P$), see
\cite[Lemma 4.4]{D}. On the other hand, $C^+(P)\otimes \Z_{\ell}$
coincides with the even Clifford $\Z_{\ell}$-algebra
$C^+(P_{\ell})$ of $(P_{\ell},\psi_{\ell})$. Clearly,
$C^+(P_{\ell})$ carries a natural $\Ga$-module structure
induced by the Galois action on $P_{\ell}$ (via the identification
$C^+(P_{\ell})=\oplus_i \wedge_{\Z_{\ell}}^{2i}P_{\ell}$). In his
adaptation of the Kuga--Satake construction, Deligne
(\cite{D}, pp. 219--223, especially Prop. 5.7 and  Lemma 6.5.1,
see also \cite{PS} and \cite{DM}, pp. 218--219) shows that 
after replacing $k$ by a finite extension, there exists 
an abelian variety $A$ over $k$,
%with complex multiplication by $C$
and an injective ring homomorphism 
$$u:C^+(P)\hookrightarrow \End(\H^1(A(\C),\Z))$$ 
satisfying the following properties.

\begin{itemize}
\item[(a)] $u:C^+(P)\hookrightarrow \End(\H^1(A(\C),\Z))$ is
a morphism of weight zero Hodge structures.

\item[(b)] The $\Z_{\ell}$-algebra homomorphism
$$u_{\ell}:C^+(P_{\ell})\hookrightarrow
\End_{\Z}(\H^1(A(\C),\Z))\otimes\Z_{\ell}=
\End_{\Z_{\ell}}(\H^1_\et(\ov A,\Z_\ell))$$ obtained from $u$ by
tensoring it with $\Z_{\ell}$, and then applying the
comparison isomorphism $\H^1(A(\C),\Z)\otimes\Z_{\ell}
=\H^1_\et(\ov A,\Z_\ell)$,
is an injective homomorphism of Galois modules.
\end{itemize}
Replacing, if necessary, $k$ by a finite extension we may
and will assume that all the endomorphisms of $\ov A$ are defined over
$k$, that is, $\End(A)=\End(\ov{A})$.

Using the compatible isomorphisms (see Subsection \ref{H1A})
$$
\begin{array}{c}
\H^1(A(\C),\Z)=\Hom(\H_1(A(\C),\Z),\Z),\\
\\
\H^1_\et(\ov A,\Z_\ell) =\Hom_{\Z_{\ell}}(T_{\ell}(A),\Z_{\ell}),
\quad  T_{\ell}(A)=\H_1(A(\C),\Z)\otimes\Z_{\ell},
\end{array}
$$
we obtain the compatible ring anti-isomorphisms
$$t:\End(\H^1(A(\C),\Z))\cong \End(\H_1(A(\C),\Z)),$$
$$t_{\ell}:\End_{\Z_{\ell}}(\H^1_\et(\ov A,\Z_\ell))\cong
\End_{\Z_{\ell}}(T_\ell(A))$$ of weight zero Hodge structures
and Galois modules, respectively. Taking the compositions, we get
an injective homomorphism of weight zero Hodge structures
$$t\, u:C^+(P)\hookrightarrow \End(\H_1(A(\C),\Z)),$$
which, extended by $\Z_{\ell}$-linearity, coincides with the
injective homomorphism of Galois modules
$$t_{\ell}\, u_{\ell}:C^+(P_{\ell})\hookrightarrow \End_{\Z_{\ell}}(T_\ell(A)).$$
We shall identify $C^+(P)$ and $\End(A)$ with their images in
$\End(\H_1(A(\C),\Z))$. Note that all the elements of
$\End (A)\subset \End(\H_1(A(\C),\Z))$ have pure Hodge
type $(0,0)$.

Let us first consider the case when $\rk\, \NS(\ov X)\geq 2$. Then
there exists a non-zero element $m\in \NS(\ov X)^\Ga\cap P$. Then
$$m\wedge T_X\subset \wedge^2 P\subset
C^+(P)\subset\End(\H_1(A(\C),\Z)).$$ Since $m\wedge T_X$ does
not contain non-zero elements of type $(0,0)$, we have
$$
(m\wedge T_X)\cap \End(A)=0.
$$
Using (\ref{K3n}) and (\ref{K3compT}), we observe that for all but
finitely many $\ell$ the $\Ga$-module $T_X/\ell$ is
isomorphic to
$$(m\wedge
T_X)/\ell\subset\End_{\Z_{\ell}}(T_{\ell}(A))/\ell=
\End_{\F_{\ell}}(A_\ell).$$
Lemma \ref{3.1} implies that $(m\wedge T_X)/\ell$
intersects trivially with $\End(A)/\ell$ for almost all
$\ell$. By the variant of the Tate conjecture 
(Proposition \ref{TLB}), for almost all $\ell$ we have
$\End_{\F_{\ell}}(A_{\ell})^\Ga=\End_{\Ga}(A_{\ell})=
\End(A)/\ell$, thus every $\Ga$-invariant element of
$m\wedge (T_X/\ell)$ is contained in
$\End(A)/\ell$, and hence must be zero. It follows that
$(T_X/\ell)^\Ga=0$ for almost all $\ell$.
Therefore, $\H^2_\et(\ov X,\mu_\ell)^\Ga=(\NS (\ov
X)/\ell)^\Ga$ for almost all $\ell$.

It remains to consider the case $\rk\, \NS(\ov X)=1$.
Then $T_X=P\simeq\Z^{21}$, and so $\wedge^{20}T_X$ is the dual
lattice of $T_X$. We have
$$\wedge^{20}T_X=
\wedge^{20} P\subset C^+(P)\subset\End(\H_1(A(\C),\Z)).$$
Since $T_X$ does not contain non-zero elements of type $(0,0)$,
the same is true for the dual Hodge structure $\wedge^{20}T_X$.
Thus $\wedge^{20}T_X\cap \End(A)=0$, and
the same arguments as before show that
$(\wedge^{20}T_X/\ell)^\Ga=0$ for almost all
$\ell$. The bilinear $\Z/\ell$-valued form induced by the
cup-product on $T_X/\ell\subset \H^2_\et(\ov X,\mu_\ell)$
is non-degenerate for almost all $\ell$, so that this
Galois module is self-dual. Thus the Galois modules
$T_X/\ell$ and $\wedge^{20}T_X/\ell$ are
isomorphic, and we conclude that $(T_X/\ell)^\Ga=0$. This
finishes the proof.
\end{proof}

\ble \label{K3a} Let $X$ be a K3 surface over a field $k$ finitely
generated over $\Q$. Then
 $\Br (\ov X)^\Ga(\ell)$ is finite for all $\ell$. \ele

\begin{proof}
By Proposition \ref{1}, it suffices to check the validity of the
Tate conjecture for divisors and the semisimplicity of the
Galois module $\H^2_\et(\ov X,\Q_\ell(1))$. Both these
assertions follow from the corresponding results on
abelian varieties, proved by Faltings in \cite{F,F86}. The
latter follows from the semisimplicity of the Galois action on
the $\ell$-adic cohomology groups of abelian varieties combined
with Proposition 6.26(d) of \cite{DM}. The former 
follows from the validity of the Tate conjecture for divisors
on abelian varieties, as explained on p. 80 of \cite{T91}.
\end{proof}

\medskip

\begin{proof}[End of proof of Theorem $\ref{t}$]
By Remark \ref{GammaIn}, it suffices to show that $\Br(\ov X)^\Ga$
is finite. By the exact sequence
(\ref{e1.1}), Corollary \ref{c} and Lemma \ref{K3} we have
$\Br(\ov X)^\Ga_\ell=0$ for almost all $\ell$. Now the
finiteness of  $\Br(\ov X)^\Ga$ follows from Lemma \ref{K3a}.
\end{proof}


\begin{thebibliography}{99}

\bibitem{Ber} V.G. Berkovich, {\em The Brauer group of abelian
varieties}.  Funktsional. Anal. i Prilozhen. {\bf 6} (1972), no. 3,
10--15; Functional Anal. Appl. {\bf 6} (1972), no. 3, 180--184.

\bibitem{CSS} J-L. Colliot-Th\'el\`ene, A.N. Skorobogatov and
Sir Peter Swinnerton-Dyer, {\em Hasse principle for pencils of
curves of genus one whose Jacobians have rational $2$-division
points}.  Inv. Math. {\bf 134} (1998), 579--650.

\bibitem{D} P. Deligne, {\em La conjecture de Weil pour les surfaces K3}.
 Inv. Math. {\bf 15} (1972), 206--226.

\bibitem{DWI} P. Deligne, {\em La conjecture de Weil}. I. Publ. Math.
IHES {\bf 43} (1974), 273--307. 

\bibitem{DW2} P. Deligne, {\em La conjecture de Weil}. II. Publ. Math.
IHES {\bf 52} (1980), 137--252.

\bibitem{DH} P. Deligne (notes by J. Milne), {\em Hodge cycles on
abelian varieties}. In: Hodge cycles, motives and Shimura
varieties. Springer Lecture Notes in Math. {\bf 900} (1982), 
9--100.

\bibitem{DM} P. Deligne, J. Milne, {\em Tannakian categories}. In:
Hodge cycles, motives and Shimura varieties. Springer Lecture
Notes in Math. {\bf 900} (1982), 101--228.

\bibitem{F} G. Faltings, {\em Endlichkeitss\"atze f\"ur abelsche
Variet\"aten \"uber Zahlk\"orpern}. Inv. Math. {\bf 73} (1983),
349--366; Erratum {\bf 75} (1984),  381.

\bibitem{F86} G. Faltings, {\em  Complements to Mordell}. In:
G. Faltings, G. W\"ustholz et al., Rational points, 2nd edition,
F. Viehweg \& Sohn, 1986.

%\bibitem{FK} E. Freitag, R. Kiehl. \'Etale cohomology and the Weil
%$conjecture. Springer-Verlag, Berlin Heidelberg New York, 1988.

\bibitem{G} O. Gabber, {\em Sur la torsion dans la cohomologie 
$\ell$-adique d'une
vari\'et\'e}.   C. R. Acad. Sci. Paris S\'er. I Math. {\bf 297}
(1983), no. 3, 179--182.

\bibitem{Gr} A. Grothendieck, {\em  Le groupe de Brauer}. I, II, III. In:
 Dix expos\'es sur la cohomologie des sch\'emas (A.
Grothendieck, N.H. Kuiper, eds), North-Holland, 1968, 46--188;
available at www.grothendieckcircle.org .


\bibitem{H} D. Harari, {\em  Obstructions de Manin transcendantes}. In:
 Number theory (Paris, 1993--1994), LMS Lecture Note Ser. {\bf
235} Cambridge Univ. Press, 1996, 75--87.

\bibitem{SGA5}J. P. Jouanolou. {\em Cohomologie de quelques sch\'emas
classiques et th\'eorie cohomologique des classes de Chern}.
Expos\'e VII. SGA 5, Cohomologie $\ell$-adique et fonctions $L$
(dirig\'e par A Grothendieck). Springer Lecture Notes in Math.
{\bf 589} (1977).

\bibitem{KleimanM} S. Kleiman, {\em The standard conjectures}.
 In: Motives (U. Jannsen, S. Kleiman, J.-P.
Serre, eds). Proc. Symp. Pure Math. {\bf 55}, Part 1 (1991),
3--20. Amer. Math. Soc., Providence, RI.

\bibitem{Kleiman} S. Kleiman, {\em The Picard scheme}. In: Fundamental
algebraic geometry (Grothendieck's FGA explained). Mathematical
Surveys and Monographs vol. {\bf 123} (2005), Amer. Math. Soc.,
Providence, RI; arXiv math.AG/0504020.

\bibitem{Lang} S. Lang, Abelian varieties, 2nd edition. Springer-Verlag, 1983.

\bibitem{Manin} Yu.I. Manin, {\em Le groupe de Brauer--Grothendieck en
g\'eom\'etrie diophantienne}. In:  Actes Congr\`es Int. Math. Nice
I (Gauthier-Villars, 1971), 401--411.

\bibitem{Milne} J.S. Milne, {\em  On a conjecture of Artin and Tate}.
Ann. Math. {\bf 102} (1975), 517--533.

\bibitem{MilnE} J.S. Milne, \'Etale cohomology. Princeton
University Press, 1980.

\bibitem{MilneAV} J.S. Milne, {\em Abelian varieties}. In:
Arithmetic geometry (G. Cornell, J.H. Silverman, eds.).
Springer-Verlag, 1986.

\bibitem{Mu} D. Mumford.  Abelian varieties, 2nd edition. Oxford University
Press, 1974.

\bibitem{PS} I.I.  Piatetski-Shapiro, I.R. Shafarevich, 
{\em  The arithmetic of surfaces
of type} K3.  Trudy Mat. Inst. Steklov. {\bf 132} (1973), 44--54;
Proc. Steklov Institute of Mathematics, {\bf 132} (1975), 45--57.

\bibitem{RS} W. Raskind and V. Scharaschkin,
{\em Descent on simply connected surfaces over algebraic number
fields}. In:  Arithmetic of higher-dimensional algebraic varieties
(Palo Alto, 2002), B. Poonen, Yu. Tschinkel, eds. Progr. Math.
{\bf 226} Birkh\"auser, 2004, 185--204.

\bibitem{SerreLie} J-P. Serre, {\em Lie Algebras and Lie Groups}.  
Springer Lecture Notes in Math. {\bf 1500} (1992). 

\bibitem{SerreLocal} J-P. Serre, {\em  Sur les groupes de Galois attach\'es aux groupes p-divisibles}. In: Proc. conf. on local fields (Driebergen, 1966),
Springer-Verlag, Berlin, 1967, 118--131.


\bibitem{S1}
A.N. Skorobogatov, {\em On the elementary obstruction to the
existence of rational points}. Mat. Zametki {\bf 81}
(2007), 112--124. (Russian)

\bibitem{SS}
A. Skorobogatov and P. Swinnerton-Dyer,  2-{\em descent on
elliptic curves and rational points on certain Kummer surfaces}.
Adv. Math. {\bf 198} (2005), 448--483.

\bibitem{SD} P. Swinnerton-Dyer, {\em Arithmetic of diagonal quartic
surfaces}, II.  Proc. London Math. Soc. {\bf 80} (2000), 513--544.

\bibitem{Tankeev} S. G.Tankeev, {\em  On the Brauer group of an 
arithmetic scheme}. II.
  Izv. Ross. Akad. Nauk Ser. Mat. {\bf 67} (2003),
155--176;    Izv. Math. {\bf 67} (2003),  1007--1029.

\bibitem{tateW} J. Tate, {\em Algebraic cycles and poles of zeta
functions}. In: Arithmetical algebraic geometry, Harper and Row,
New York, 1965, 93--110.

\bibitem{tateB} J. Tate, {\em  On the conjectures of Birch and
Swinnerton-Dyer and a geometric analog}, S\'eminaire Bourbaki
{\bf 306} (1965/1966); In: Dix expos\'es sur la cohomologie des
sch\'emas (A. Grothendieck, N.H. Kuiper, eds), North-Holland,
1968, 189--214.

\bibitem{tate} J. Tate, {\em Endomorphisms of abelian varieties over finite
fields}.  Inv. Math. {\bf 2} (1966), 134--144.

\bibitem{Tate76} J. Tate, {\em  Relations between $K_2$ and Galois
cohomology}.  Inv. Math. {\bf 36} (1976), 257--274.

\bibitem{T91} J. Tate, {\em  Conjectures on algebraic cycles in
$\ell$-adic cohomology}. In: Motives (U. Jannsen, S. Kleiman,
J.-P. Serre, eds). Proc. Symp. Pure Math. {\bf 55}, Part 1 (1991),
71--83. Amer. Math. Soc., Providence, RI.

\bibitem{We}R.O. Wells, Differential analysis on complex
manifolds. Prentice-Hall, Inc., Englwood-Cliffs, NJ, 1973.

\bibitem{W} O. Wittenberg, {\em  Transcendental Brauer--Manin obstruction
on a pencil of elliptic curves}. In:  Arithmetic of
higher-dimensional algebraic varieties (Palo Alto, 2002), B.
Poonen, Yu. Tschinkel, eds. Progr. Math. {\bf 226} Birkh\"auser,
2004, 259--267.


\bibitem{Z75} Yu.G. Zarhin, {\em  Endomorphisms of Abelian varieties 
over fields of
finite characteristic}. Izv. Akad. Nauk SSSR Ser. Mat. 
{\bf 39} (1975), 272--277; Math.
USSR Izv. {\bf 9} (1975), 255--260.

\bibitem{Z76} Yu.G. Zarhin, {\em Abelian varieties in characteristic} $p$.
 Mat. Zametki {\bf 19} (1976), 393--400; Math. Notes {\bf 19} (1976), 240--244.

\bibitem{Z77} Yu.G. Zarhin, {\em  Endomorphisms of abelian varieties
and points of finite order in characteristic} $p$. Mat. Zametki
{\bf 21}  (1977), 737--744; Math. Notes {\bf 21} (1977), 415--419.

\bibitem{Z82} Yu.G. Zarhin, {\em  The Brauer group of an abelian variety
over a finite field}.  Izv. Akad. Nauk SSSR Ser. Mat. {\bf 46}
(1982), 211--243; Math. USSR Izvestia {\bf 20} (1983), 203--234.

\bibitem{Z85} Yu.G. Zarhin, {\em  A finiteness theorem
for unpolarized abelian varieties over number fields with
prescribed places of bad reduction}.  Inv. Math. {\bf 79} (1985),
309--321.

\end{thebibliography}
\end{document}